\documentclass[twocolumn,secnumarabic,amssymb, nobibnotes, aps, prd]{revtex4-2}

\usepackage{amscd}
\usepackage{amsmath}
\usepackage{amssymb}
\usepackage{graphicx}%
\usepackage{dsfont}
\usepackage{bm}

\usepackage{float}
\usepackage{placeins}

\usepackage{newtxtext}
\usepackage{booktabs,caption} 
\usepackage{sidecap}
\usepackage{grffile}

\usepackage[colorlinks=true, pdfstartview=FitV, linkcolor=blue,
            citecolor=blue, urlcolor=blue]{hyperref}

\usepackage{xcolor}

\usepackage{scalerel}

\graphicspath{{./Figures/}}

\def\bi{\begin{itemize}}
\def\ei{\end{itemize}}
\def\be{\begin{equation}}
\def\ee{\end{equation}}
\def\bea{\begin{equation} \begin{aligned}}
\def\eea{\end{aligned} \end{equation}}
\def\beas{\begin{equation*} \begin{aligned}}
\def\eeas{\end{aligned} \end{equation*}}
\def\bes{\begin{equation*}}
\def\ees{\end{equation*}}

\def\d{\, \mathrm{d}}

\def\x{\boldsymbol{x}}
\def\y{\boldsymbol{y}}

\def\z{\boldsymbol{z}}

\newcommand{\T}{\boldsymbol{\theta}}
\def\X{\boldsymbol{\mathcal{X}}}

\def\V{\boldsymbol{\mathcal{V}}}
\def\Z{\boldsymbol{\mathcal{Z}}}
\def\cN{\mathcal{N}}

\newcommand{\norm}[1]{\left\lVert#1\right\rVert}

\begin{document}

\title{The High-Frequency and Rare Events Barriers to Neural Closures of Atmospheric Dynamics}%

\author{Micka\"el D. Chekroun}
\email{mchekroun@atmos.ucla.edu} 
\affiliation{Department of Atmospheric and Oceanic Sciences, University of California, Los Angeles, CA}
\affiliation{Department of Earth and Planetary Sciences, Weizmann Institute of Science, Rehovot 76100, Israel}

\author{Honghu Liu}
\affiliation{Department of Mathematics, Virginia Tech, Blacksburg, VA 24061, USA}

\author{Kaushik Srinivasan}
\affiliation{Department of Atmospheric and Oceanic Sciences, University of California, Los Angeles, CA 90095-1565, USA}

\author{James C. McWilliams}
\affiliation{Department of Atmospheric and Oceanic Sciences and Institute of Geophysics and Planetary Physics, University of California, Los Angeles, CA 90095-1565, USA} 

\begin{abstract}
Recent years have seen a surge in interest for leveraging neural networks to parameterize small-scale or fast processes in climate and turbulence models.
In this short paper, we point out two fundamental issues in this endeavor. 
The first concerns the difficulties neural networks may experience in capturing rare events due to limitations in how data is sampled. The second arises from the inherent multiscale nature of these systems. They combine high-frequency components (like inertia-gravity waves) with slower, evolving processes (geostrophic motion). This multiscale nature creates a significant hurdle for neural network closures.
To illustrate these challenges, we focus on the atmospheric 1980 Lorenz model, a simplified version of the Primitive Equations that drive climate models. This model serves as a compelling example because it captures the essence of these difficulties.
  \end{abstract}

\date{March 18, 2024}%
\maketitle

\section{Introduction}
Atmospheric and oceanic flows constrained by Earth's rotation
   satisfy an approximately geostrophic momentum balance on larger
   scales, associated with slow evolution on time scales of days, but
   they also exhibit fast inertia-gravity wave oscillations.  The
   problems of identifying the slow component (e.g., for weather
   forecast initialization \cite{bolin1955numerical,baer1977complete,machenhauer1977dynamics,daley1981normal}) and of characterizing slow-fast
   interactions are central to geophysical fluid dynamics, and the former was first coined as a slow manifold problem by Leith \cite{leith1980nonlinear}.  
   The L63 model \cite{lorenz1963deterministic} famous for its chaotic strange attractor is a paradigm for the
   geostrophic component, while the L80 model \cite{Lorenz80}  is its paradigmatic successor both
   for the generalization of slow balance and for slow-fast coupling.

The explosion of machine learning (ML) methods provides an unprecedented opportunity to analyze data and accelerate scientific progress. A variety of ML methods have emerged for solving dynamical systems \cite{sirignano2018dgm,raissi2019physics,bar2019learning}, predicting \cite{pathak2018model} or discovering \cite{brunton2016discovering} them from data. For larger scale problems, much effort has been devoted lately to the learning of neural subgrid-scale parameterizations in coarse-resolution climate models \cite{rasp2018deep} but yet the lack of interpretability and reliability prevents a widespread adoption so far \cite{gentine2018could,brenowitz2020interpreting}.
   
In parallel, the learning of stable neural parameterizations of small scales or neglected variables has progressed remarkably for the closure of fluid models in turbulent regimes such as the forced Navier-Stokes equations or quasi-geostrophic flow models on a $\beta$-plane; see \cite{bolton2019applications,maulik2019,kochkov2021,zanna2020,subel2022explaining,CNN_small23,Lucarini_Chekroun2023}. 

While neural networks show promise for climate modeling, the full Primitive Equations (PE) remain a challenge. This study identifies potential hurdles in achieving efficient neural closures for PE. We leverage the L80 model, a simplified version of the PE, as a illustrative example to highlight these fundamental issues.

 In that respect, the L80 model exhibits a fascinating dynamical transition. For small Rossby numbers, its solutions evolve slowly over time and remain entirely slow, dominated by large-scale Rossby waves \cite{Gent_McWilliams82}. However, as the Rossby number increases, faster oscillations become superimposed on these slow background motions \cite{vautard1986invariant,CLM16_Lorenz9D}. This spontaneous emergence of high-frequency components, linked to inertia-gravity waves (IGWs) riding on the slower geostrophic flow, significantly complicates the closure problem in atmospheric models \cite{CLM16_Lorenz9D,CLM21_BE}.

Multiscale dynamics, characterized by the intricate interplay of slow and fast processes without clear separation, are not unique to the L80 model. Similar regimes have been observed in fully resolved Primitive equation (PE) models, where fronts and jets generate complex multiscale interactions \cite{plougonven2007inertia,polichtchouk2020spontaneous} as well as in cloud-resolving models, where large-scale convectively coupled gravity waves emerge spontaneously \cite{tulich2007vertical}. Tropical convection regions, where organized activity produces gravity waves with a broad spectrum, ranging from 10 km to over 1000 km wavelengths \cite{lane2015gravity} provide another instance of such multiscale dynamics. Finally, inertia-gravity waves have also been observed in continental shallow convection, where they contribute to organized mesoscale patterns over vegetated areas  \cite{Dror2021}.

Inertia-gravity waves  can hold surprising amounts of energy even at large scales. 
For example, Rocha et al.~\cite{rocha2016mesoscale} found that IGWs contribute nearly half of the near-surface kinetic energy in specific ocean regions at scales ranging from 10 to 40 km. This overlap between wave and turbulence scales in geophysical kinetic energy spectra creates a challenge: perturbation methods like Wentzel-Kramers-Brillouin (WKB) \cite{bender1999advanced} become inapplicable across all scales \cite{young2021inertia}.

Such regimes where slow and fast dynamics overlap were shown to constitute critical challenges for closure methods in the L80 model. Solutions in these regimes blend slow background motion with sudden bursts of IGWs carrying a significant portion of the total energy. These ``high-low frequency (HLF)" solutions disrupt the expected slaving relationships satisfied at lower Rossby numbers, leading to a major breakdown in closure techniques relying on a separation between the slow and fast variables \cite{CLM16_Lorenz9D}.

A recent study by \cite{CLM21_BE} proposes a promising solution to closure problems in such HLF regimes without timescale separation and where slow Rossby variables are influenced by high-frequency waves. This approach hinges on the Balance Equation (BE) \cite{mcwilliams1980intermediate,Gent_McWilliams82}  as rooted in the works of Monin \cite{Monin1952}, Charney and Bolin \cite{charney1955use,bolin1955numerical}, and Lorenz \cite{lorenz1960energy}, which allows for a nonlinear separation of variables. As demonstrated in \cite{CLM21_BE}, the BE isolates, for large Rossby numbers, the fast, non-geostrophic component of the flow as residual dynamics off the BE manifold. Building on the BE separation, it was shown in  \cite{CLM21_BE} that this fast motion can be effectively parameterized using networks of nonlinear stochastic oscillators (NSOs).  These NSOs are designed to match the characteristic patterns of variability observed in the fast motion, leveraging the concept of resonances discussed in  \cite{Chek_al14_RP,Chekroun_al_RP2,RP_Hopf}.  
The resulting stochastic closure shows then high-accuracy skills in reproducing the multiscale dynamics.

This work emphasizes the limitations of (standard) neural networks (alone) for achieving such accurate closures for HLF regimes,  highlighting their struggle  to simultaneously capture the slow, balanced motion while restoring the high-frequency oscillations.
Section \ref{Sec_sensitivity} discusses the limitations of neural networks for parameterizing the L80 model's slow motion,  emphasizing in particular their sensitivity to rare event statistics (Section \ref{Sec_slow_erg}). Section \ref{Sec_highfreq} highlights the fundamental challenges faced by neural networks in capturing both the slow and high-frequency content of the L80 solutions, ultimately hindering accurate closure.

\section{Learning slow neural closure: Sensitivity}\label{Sec_sensitivity}

The L80 model, obtained by Lorenz  in \cite{Lorenz80} as a nine-dimensional truncation of the PE onto three Fourier modes with low wavenumbers, can be written as:
\bea\label{Eq_L9D}
 \!\! a_i \frac{\d x_i}{\d t} & = -\nu_0 a_i^2 x_i - c(a_i - a_k) x_j y_k + c(a_i - a_j) y_j x_k\\
 & \quad + a_i b_ix_jx_k  - 2c^2y_jy_k  + a_i(y_i - z_i), \\
 \!\! a_i \frac{\d y_i}{\d t} & =  -  a_kb_k x_jy_k - a_jb_j y_jx_k + c(a_k-a_j)y_jy_k\\
 & \quad -a_ix_i-\nu_0a_i^2y_i, \\
 \!\! \frac{\d z_i}{\d t} & = g_0 a_ix_i - b_kx_j(z_k-h_k) - b_j(z_j-h_j)x_k\\
 & + cy_j(z_k-h_k) - c(z_j-h_j)y_k -\kappa_0a_iz_i +  F_i,
 \eea
whose model parameters are described in \cite{Lorenz80,CLM16_Lorenz9D}.
 
The above equations are written for each cyclic permutation of the set of indices $(1, 2, 3)$, namely, for  $(i, j, k)$ in $\{(1, 2, 3), (2, 3, 1), (3, 1, 2)\}$.  The model variables $(\x,\y,\z)$ are amplitudes for the divergent velocity potential, stream-function, and dynamic height, respectively.

In this model, the square root of the constant forcing $F_1$ can be interpreted as the Rossby number; see \cite{Gent_McWilliams82} and \cite[Eq.~(2.4)]{CLM16_Lorenz9D}. Transitions to chaos occur as the Rossby number $Ro$ is increased \cite{Gent_McWilliams82,CLM16_Lorenz9D}. As mentioned above, at small Rossby numbers, the solutions to the L80 model are dominated by Rossby waves and thus remain entirely slow for all time. As identified in \cite{CLM16_Lorenz9D}, when the Rossby number is further increased beyond a critical Rossby number $Ro^\ast$,  fast IGW oscillations emerge spontaneously and are superimposed on the slow component of the solutions. For such regimes, the aforementioned BE manifold on which balanced solutions lie \cite{mcwilliams1980intermediate,Gent_McWilliams82,CLM16_Lorenz9D} is no longer able to parameterize fully the L80 dynamics since a substantial portion of it, associated with the IGWs, evolves transversally to the BE manifold \cite[Fig.~3]{CLM21_BE}. These regimes with energetic bursts of IGWs lie beyond the parameter range explored by Lorenz in his original 1980 article \cite{Lorenz80} and beyond other regimes with exponential smallness of IGW amplitudes as studied in subsequent Lorenz 86 models \cite{lorenz1986existence,lorenz1987nonexistence,camassa1995geometry,vanneste2008exponential} and the full primitive equations \cite{temam2011slow} at smaller Rossby numbers \cite{vanneste2013balance}.

\begin{figure*}[htb!]
\centering
\includegraphics[width=0.75\textwidth,height=.55\textwidth]{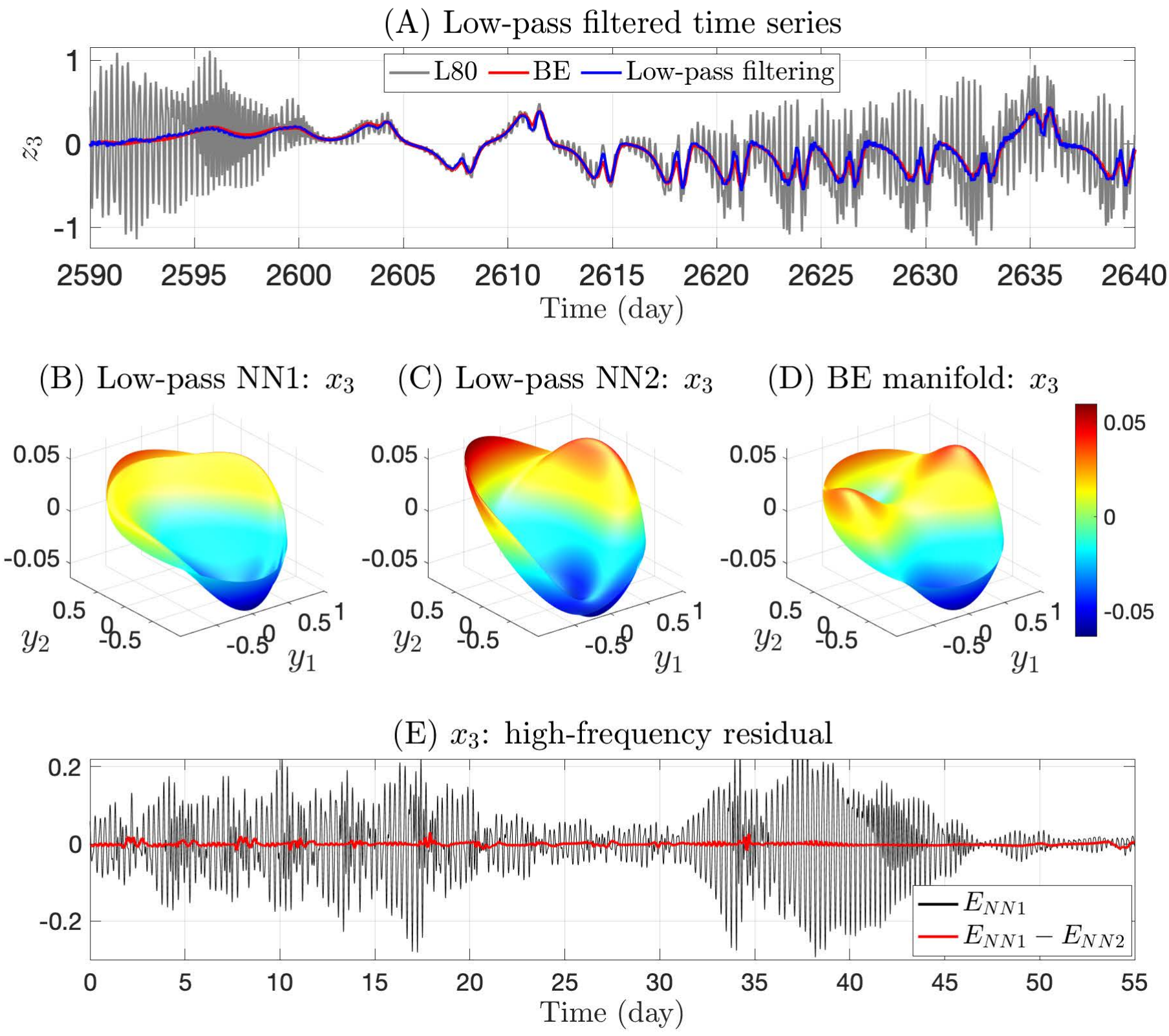}
\caption{{\bf Panel A}: Illustration, for the $z_3$-variable, of the BE manifold's ability in capturing the L80 model's slow motion. See \cite{Gent_McWilliams82} and Appendix \ref{Sect_Appendix_BE} for a derivation.  
{\bf Panels (B) and (C):}  Neural parameterizations  $\X_{\T}^3$ for the $x_3$-variable, as learnt through random selection (NN$_1$)/predefined selection (NN$_2$). Visualized here as mappings from $(y_1,y_2)$ onto the unit sphere in $\mathbb{R}^3$. {\bf Panels (D):} Same visualization adopted for the BE manifold.
 {\bf Panels (E)}: High-frequency residual $E_{NN_1}(t)$ for $x_3$ (black) given by \eqref{Eq_residual_err} and its difference with $E_{NN_2}(t)$ (red).  
}\label{Fig_sensitivity}
\end{figure*}

 The HLF solutions considered in this study are obtained for such a critical parameter regime where $Ro >Ro^\ast$. They correspond to those of \cite[Fig.~7]{CLM21_BE}; see Appendix~\ref{Sect_Appendix_HLF} for details. We first analyze the ability of neural parameterizations to learn the slow motion of the L80 dynamics in the HLF regime.
To do so, we preprocess the target variables  $\x$ and $\z$ to be parameterized
by applying a low-pass filter in order to extract the slow motion.     
In that respect, a simple moving average is adopted with a window size equal to $T_{GW}$, the dominant period of the gravity waves.
The results are shown in Figure \ref{Fig_sensitivity}A for the $z_3$-variable for which we observe that the low-pass filtered solution almost coincides this way with the BE parameterization $\z_{BE}(t)=G(\y(t))$ with  $\y(t)$ denoting the $\y$-component of the  HLF solution to the L80 model. 

The L80 model has an inherent structure that can be exploited for closure. Studies have shown that the BE manifold, constructed in two steps (parameterizing $\z$ as a function of $\y$ and then $\x$ as a function of $\y$ and the parameterized $\z$), achieves excellent closure across various parameter regimes \cite{Gent_McWilliams82} (see Appendix \ref{Sect_Appendix_BE} and \cite{CLM16_Lorenz9D} for details). To leverage this existing knowledge and facilitate comparison with the BE manifold, we design our neural network parameterizations with a similar structure.
Specifically, we first learn a feedforward neural network (multilayer perceptron, MLP) denoted as $\Z_{\T}$, which takes the (unfiltered) variable $\y$ as input and predicts the filtered $\z$-variable (Eq.~\eqref{Eq_Z}). Then, we train a second MLP, $\X_{\T}$, that takes both $\y$ and the output of $\Z_{\T}$, $(\y, \Z_{\T}(\y))$, as input in order to predict the filtered $\x$-variable.

The structure of our MLPs is standard.
Each neural parameterization, e.g.~$\z$ in terms of $\y$, is sought by means of an MLP with $L$ hidden layers of $p$ neurons each. It boils down to find
\be\label{Eq_Z}
\Z_{\T}(\y) = {\cN}_{\rm out} \circ {\cN}_L \circ \cdots \circ {\cN}_1 \circ {\cN}_{\rm in}(\y),
\ee
in which ${\cN}_{\rm in}$ (resp.~${\cN}_{\rm out}$) constitutes the input (resp.~output) layer, 
while ${\cal N}_k$ is a mapping from $\mathbb{R}^p$ (the space of neurons) onto itself, given by 
${\cal N}_k(\xi ) = \Psi_k({\bf W}_k \xi + {\bf b}_k)$ ($\xi$ in $\mathbb{R}^p$) where $\Psi_k$ is  
a $p$-dimensional elementwise function, i.e.~a function that applies a (scalar) activation function to each of its inputs individually,  
and the ${\bf W}_k$ and $ {\bf b}_k$ denote respectively the weight matrices and bias vectors to be learnt. In \eqref{Eq_Z}, the subscript $\T$ denotes the collection of these parameters.
In this work, the nonlinear activation function is a simple $\tanh$ function, and the input and output layers consist just of linear normalization and reversal operations. It turns out that NNs with one hidden layer and 5 neurons are sufficient to obtain loss functions with a small residual; see Table~\ref{tab_loss_func}.

\begin{table*}[htb!] 
\caption{{\bf Loss function evaluations for two neural networks}. The loss functions \eqref{Eq_loss1} for $\z$ and \eqref{Eq_loss2} for $\x$, are minimized using two neural networks,  NN$_1$ and NN$_2$ providing  each a parameterization  $(\Z_{\T},\X_{\T})$,  differing only in the way the training, validation, and testing sets are selected. In each case, the aspect ratios between these sets are the same.}\label{tab_loss_func}  
\begin{tabular}{c| c c c c c c  c}
\toprule
Epochs & 10 & 50 & 100 & 300 & 500 & 1000  \\
\hline
 NN$_1$ loss for $\z$ (random) ($\times 10^{-3}$)  & $11.17$ & $9.26$ & $9.26$ & $9.26$ & $9.26$ & 9.26\\
 NN$_2$ loss for $\z$ (predefined) ($\times 10^{-3}$) & $13.70$ & $10.66$ & $9.28$ & $9.05$ & $9.05$ & 9.05\\
NN$_1$ loss for $\x$ (random) ($\times 10^{-4}$)  & $1.76$ & $1.38$ & $1.35$ & $1.33$ & $1.32$ & 1.32\\
NN$_2$ loss  for $\x$ (predefined) ($\times 10^{-4}$) & $1.62$ & $1.37$ & $1.33$ & $1.31$ & $1.31$ & 1.31\\
\bottomrule 
\end{tabular}
\end{table*}
Based on our approach paralleling the BE manifold construction, we learn our neural parameterizations for the L80 model, through the following consecutive minimizations. First, given a discrete set of time instants $t_j$, one minimizes
\bea\label{Eq_loss1}
\mathcal{L}_{\T}(\z;\y)=\sum_j \norm{\z_{t_j}-\Z_{\T}(\y_{t_j})}^2,
\eea
in which  $\z$ is filtered (in time) while  $\y$ is not, followed by the minimization of
\bea\label{Eq_loss2}
\hspace{-.5ex}\mathcal{L}_{\T}(\x;(\y,\Z_{\T_1^\ast}(\y)))\hspace{-.3ex}=\hspace{-.3ex}\sum_j \norm{\x_{t_j}\hspace{-.6ex}-\hspace{-.4ex}\X_{\T}(\y_{t_j},\Z_{\T_1^\ast}(\y_{t_j}))}^2,
\eea
with $\x$ filtered and where $\Z_{\T_1^\ast}$ denotes the optimal parameterization obtained after minimization of \eqref{Eq_loss1}. 

\begin{figure}[hbt!]
\includegraphics[width=.49\textwidth,height=.28\textwidth]{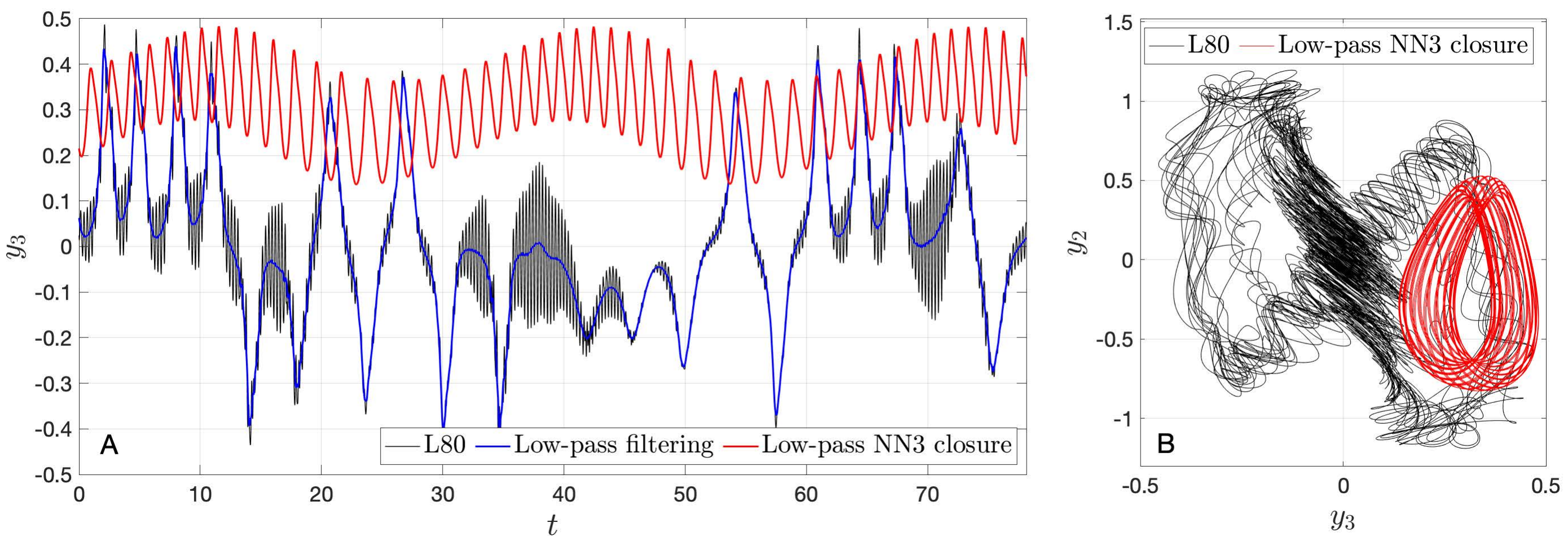}
\caption{{\bf False quasiperiodicity produced by a slow neural closure}. Here, the slow neural closure Eq.~\eqref{Eq_low-passNN} is driven by $\Z_{\T}$ and $\X_{\T}$ that are trained using a low-pass filtered version of $\y(t)$ (blue curve in Panel (A)) unlike the closure defined in Eq.~\eqref{Eq_low-passNN} where the slow neural closure is trained using the unfiltered $\y$-variable.} \label{Fig_y3_timeseries}
\end{figure}

We emphasize the importance  of including the unfiltered $\y$-component of the HLF solution in the training data, even though it contains rapid oscillations. This unfiltered data is indeed crucial for the network to learn a proper representation of the slow motion.  If we replace the unfiltered $\y$-component with a filtered version (like the blue curve for $y_3$ in Figure \ref{Fig_y3_timeseries}A), the resulting closure fails. It produces an unrealistic quasi-periodic behavior that does not resemble even the L80 model's quasi-periodic behaviors documented in \cite{Gent_McWilliams82} for nearby parameter settings (see red curves in Figure \ref{Fig_y3_timeseries}).

To assess whether a neural parameterization is successful in capturing the slow motion, we evaluate also the following {\it high-frequency (HF) residual}
\be\label{Eq_residual_err} 
E_{NN}^j(t)=x_j(t)-\X_{\T_2^\ast}^j(\y(t), \Z_{\T_1^\ast}(\y(t))),
\ee
in which the $x_j(t)$ and $\y(t)$ are {\it both} unfiltered.  For an NN with small residual, $E_{NN}^j(t)$  is typically void of slow oscillations (see Figure \ref{Fig_sensitivity}E) with mean $\langle E_{NN}^j \rangle \approx 0$ for each $1\leq j\leq 3$.  

Figure \ref{Fig_sensitivity}  illustrates this feature with two neural networks, NN$_1$ and NN$_2$, trained using different strategies for selecting training, validation, and testing data. Even though both networks achieve good parameterization results offline (similar to the BE manifold), their underlying structures differ visually from the BE manifold.

 To explore these differences, we focus on specific components ($\X_{\T_2^\ast}^j$ for $x_j$ and $\Z_{\T_1^\ast}^j$ for $z_j$) of the neural parameterizations. We plot these components as level sets on a three-dimensional sphere to reveal their geometric properties. This visualization is particularly useful since $\Z_{\T_1^\ast}^j$ and $\X_{\T_2^\ast}^j$ are scalar fields depending on three variables. For a given radius, the level sets of $\Z_{\T_1^\ast}^j$ (resp.~$\X_{\T_1^\ast}^j$) on the three-dimensional  sphere, $y_1^2+y_2^2+y_3^2=r^2$, can be visualized as a 2D surface that maps $(y_1,y_2)$ to  $z_j$ (resp.~$x_j$).  Figures \ref{Fig_sensitivity}B, \ref{Fig_sensitivity}C, and \ref{Fig_sensitivity}D show these level sets for radius $r=1$. 
 
 Interestingly, these visualizations reveal significant differences in the minimizers (and consequently, the parameterization formulas) of NN$_1$ and NN$_2$, even though their loss function values differ only by 1\% (Table \ref{tab_loss_func}) and their high-frequency residuals are similar (red curve in Figure \ref{Fig_sensitivity}E).

These geometric offline differences hide more profound consequences when the neural parameterizations are used online, for closure. As explained below, the sensitivity of online predictions that are tied to sampling issues is indeed observed.  In that respect, recall that a common practice to train NNs is to divide the dataset into three subsets. The first subset is the training set, which is used for computing the loss function's gradient and updating the network weights and biases.

The second subset is the validation set. It corresponds to the second dataset over which the prediction skills of the fitted model are assessed.  
The error on the validation set is monitored during the training process to provide an unbiased evaluation while tuning the model's hyperparameters.
When the network begins to overfit the data, the error on the validation set typically begins to rise after an initial decrease. The network parameters are saved at the minimum.  It gives then the ``final model" that is tested over the test set that is typically a holdout dataset not used as a validation nor a training set.

The  parameterization NN$_1$ shown in Figure \ref{Fig_sensitivity}A  is learnt through a random selection while NN$_2$ is learnt through a predefined selection. In each case,  ratios for training, testing, and validation are 0.7, 0.15, and 0.15, respectively.
The total length of the training is 700 days. Given the same input and target data, the minimal values of the loss functions \eqref{Eq_loss1}-\eqref{Eq_loss2} for NN$_1$ and NN$_2$ are reported in Table~\ref{tab_loss_func}, across epochs. Already after 500 epochs, one observes that the loss function evaluations differ only by $1\%$ between the random or predefined selection protocol of the training, validation, and testing sets.

\begin{figure}[ht]
\includegraphics[width=0.5\textwidth,height=.25\textwidth]{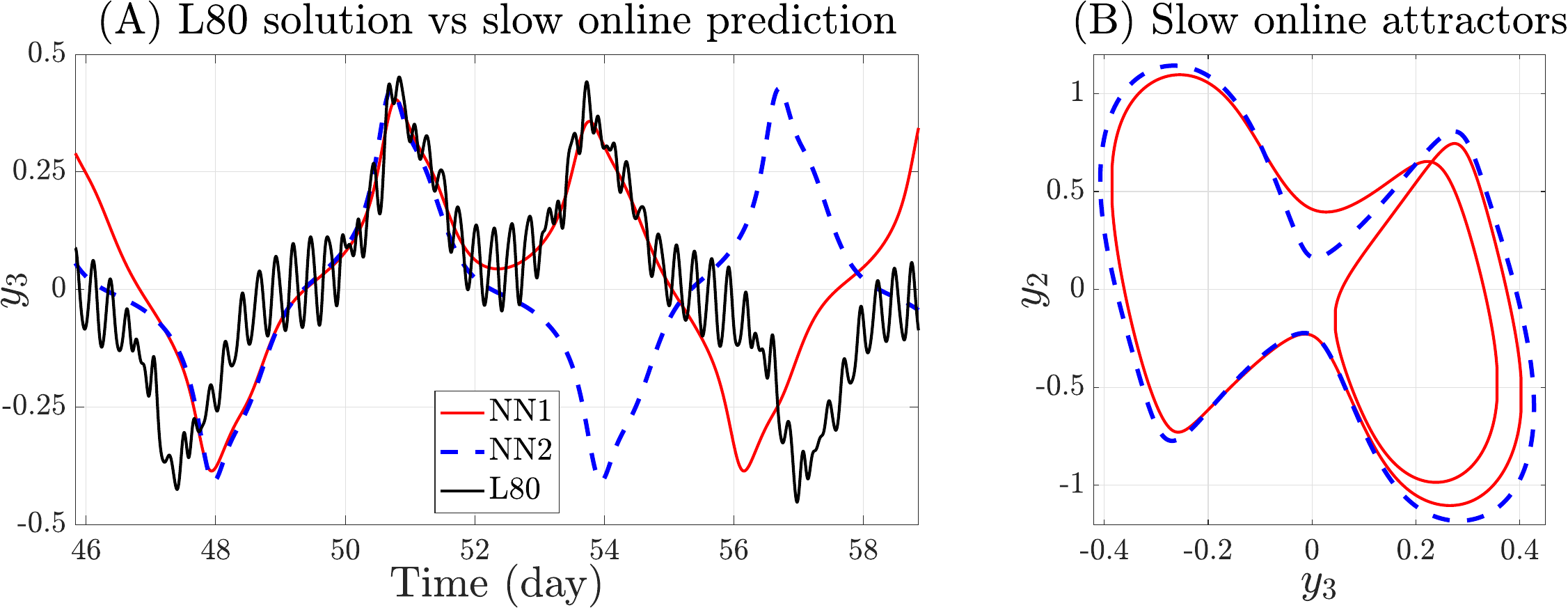} 
\caption{{\bf Sensitivity of the slow neural closures.}  Here, NN$_1$ and NN$_2$ differ only in their training modalities. NN$_1$ is learnt from random selection of the training, validation, and testing sets, and NN$_2$ from a predefined selection with the same aspect ratios; see Text. The corresponding loss functions differ by 1$\%$ (see Table~\ref{tab_loss_func}), while the dynamical differences of the online predictions are substantial.}\label{Fig_BE_vs_slowNN_closures}
\end{figure}

 We now discuss the sensitivity issue of online predictions driven by such neural parameterizations that are close in terms of their loss function scoring.
This point is illustrated in Figure \ref{Fig_BE_vs_slowNN_closures}. There, we show online prediction corresponding to a given  slow NN-parameterization $(\X_{\T_2^\ast},\Z_{\T_1^\ast})$  learnt by minimization of the loss functions (Eqns.~\eqref{Eq_loss1} and \eqref{Eq_loss2}), namely the solution to the slow neural closure
\begin{widetext}
\bea\label{Eq_low-passNN}
 a_i\frac{\d y_i}{\d t} &=  -  a_kb_k \X_{\T_2^\ast}^j(\y,  \Z_{\T_1^\ast}(\y)) y_k - a_j b_j y_j   \X_{\T_2^\ast}^k(\y,  \Z_{\T_1^\ast}(\y))+ c(a_k-a_j)y_j y_k -  a_i  \X_{\T_2^\ast}^i(\y,  \Z_{\T_1^\ast}(\y))- \nu_0a_i^2 y_i.
 \eea
 \end{widetext}
This closed equation in the $\y$-variable is obtained by replacing the $x_\ell$-variables  in the $\y$-equation of  the L80 model  (Eq.~\eqref{Eq_L9D}) by their neural parameterizations, either NN$_1$ or  NN$_2$.

The attractor corresponding to the slow NN$_1$-closure (with random selection) differs clearly from that of slow NN$_2$-closure (with predefined selection) in spite of convergence and closeness of the loss functions at their respective minimal value; see Figure \ref{Fig_BE_vs_slowNN_closures}B. Both predict periodic orbits with different attributes, one self-intersecting in the $(y_2,y_3)$-plane (NN$_1$), the other without intersection point (NN$_2$). 

A closer inspection at these topological differences reveals
in the time domain that the   
slow NN$_1$-closure is able to capture more accurately the low-frequency content of certain temporal patterns exhibited by the  HLF  solutions of the L80 model compared to the slow NN$_2$-closure; blue vs red curves in Figure \ref{Fig_BE_vs_slowNN_closures}A. We argue below that such a sensitivity between online solutions takes its root in the rare events tied to the irregular transitions exhibited by the  HLF solutions to the L80 model that spoils the offline learning.

\FloatBarrier 

\begin{figure*}[bth!]
\includegraphics[width=1\textwidth,height=.5\textwidth]{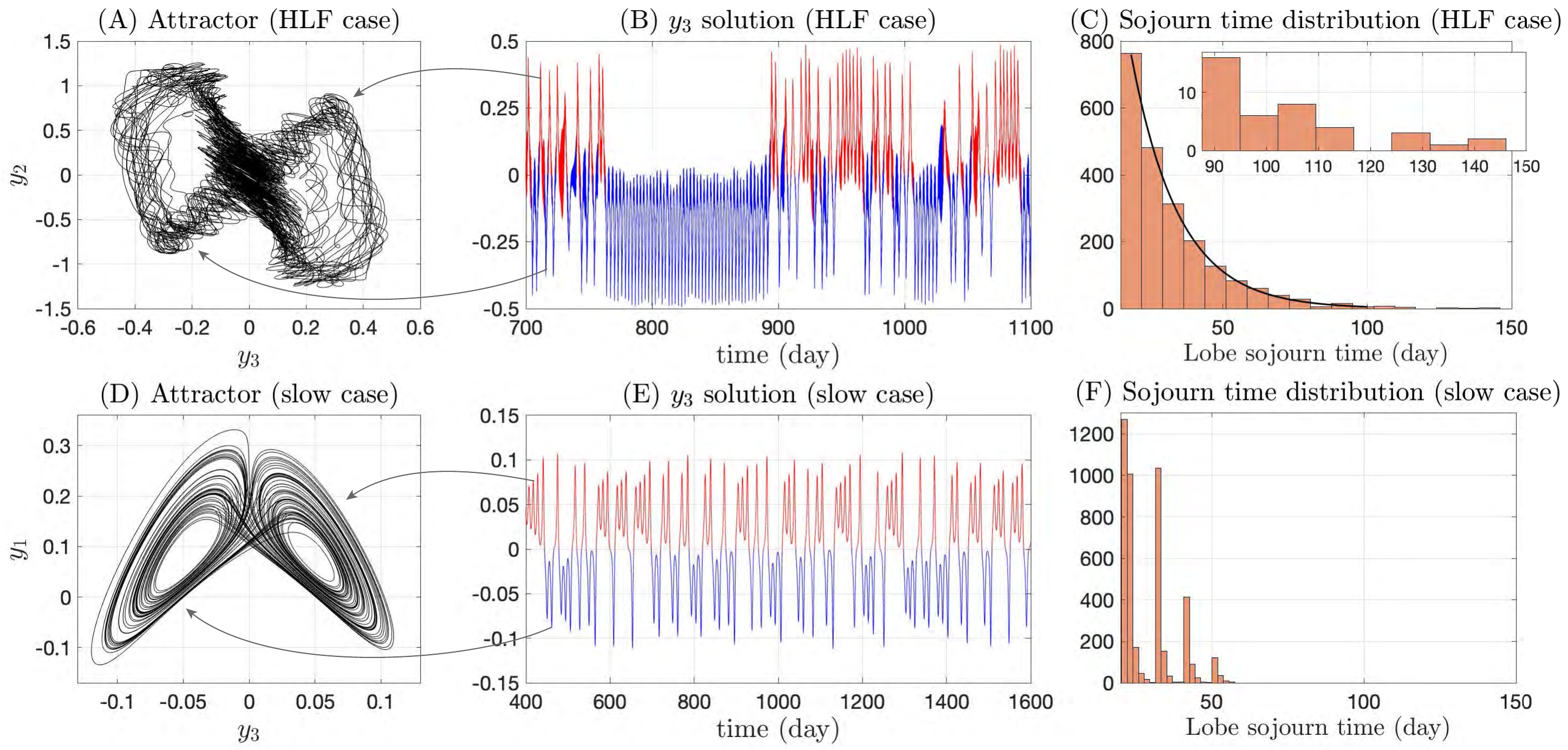}
\caption{{\bf Panel A:} Attractor in the HLF case.  {\bf Panel B:} The sojourn episodes within one particular lobe are marked by different colours. Here, the parameters are those used in Lorenz's original paper \cite{Lorenz80} except  $F_1=0.3027$  in Eq.~\eqref{Eq_L9D}. {\bf Panel C: Lobe sojourn time distributions.} The exponential fit is calculated over 500 yr-long simulation of  Eq.~\eqref{Eq_L9D} and is shown by the black curve $f(t)=a e^{b t}$ with $a=2292$ and $b=-6.05\times 10^{-2}$ with $t$ in day. The inset in panel C shows a magnification of the distribution for the rare and large sojourn times.  {\bf Panels E and F:} Same as panels B and C except that $F_1=6.97\times 10^{-2}$, corresponding to the slow chaotic regime shown in panel D in which the solutions are void of fast oscillations. In this regime, no rare event statistics emerge.}\label{Fig_sojourn}
\end{figure*}

In contrast, at lower Rossby numbers, for regimes devoid of fast oscillations such as shown  in Figure \ref{Fig_sojourn}D below  corresponding to $F_1=6.97\times 10^{-2}$ in the L80 model, neural closures of high-accuracy are easily accessible with skills comparable to those obtained with the BE manifold; see Figure \ref{Fig_slowL80_vs_NN}.  As explained below, the reasons for this success lie in the absence of high-frequencies in the solutions to parameterize and in the absence of rare events in the statistics of lobe transitions.

\section{Irregular transitions, rare events and learning consequences}\label{Sec_slow_erg}
\begin{figure}[tbh!]
\includegraphics[width=0.48\textwidth,height=.25\textwidth]{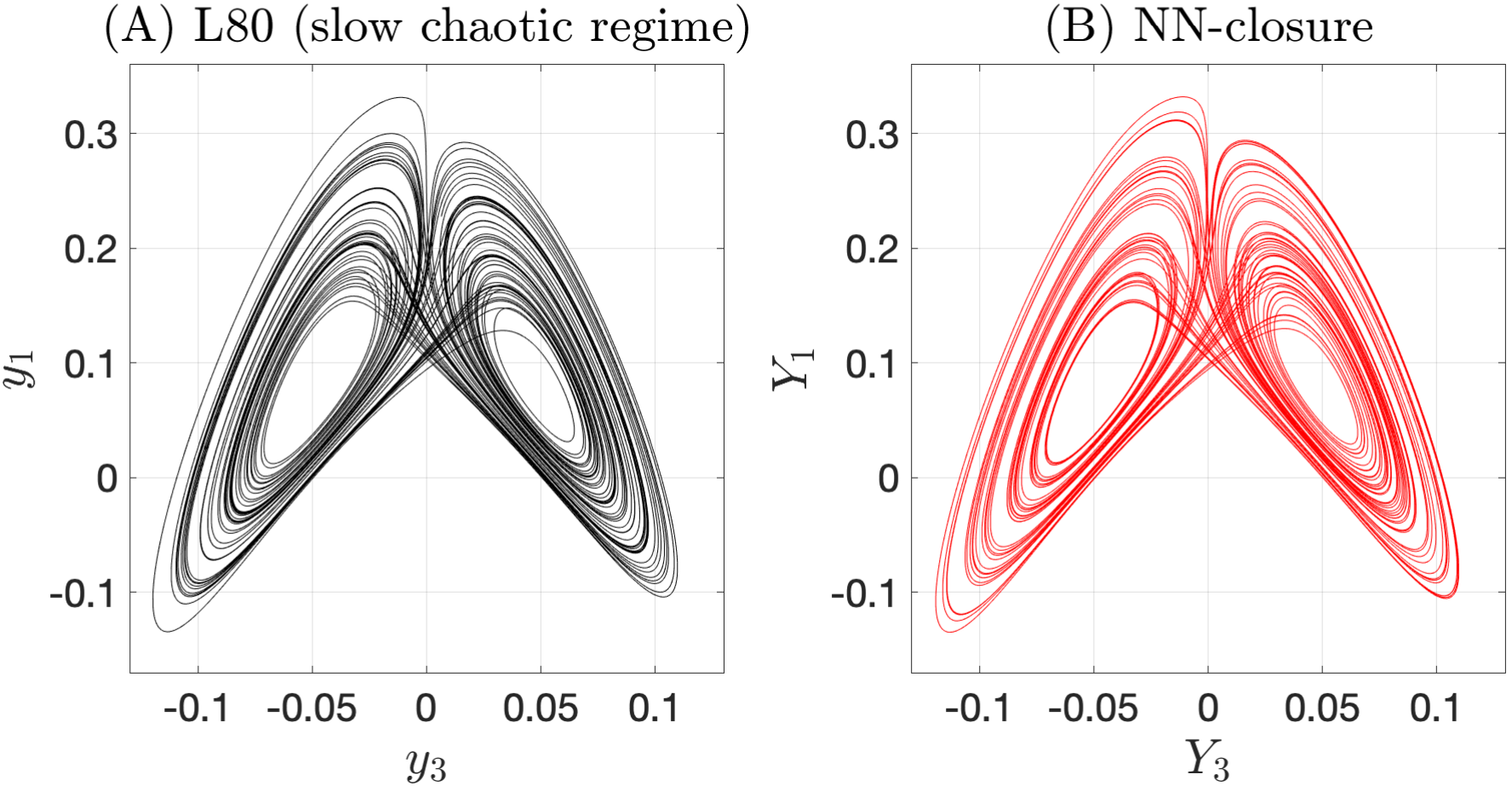}
\caption{{\bf The L80 attractor vs. its NN-closure in the slow chaos regime.}  Here $F_1=6.97\times 10^{-2}$ in the L80 model, which corresponds to the slow chaos case shown in Figure \ref{Fig_sojourn}D and in \cite[Fig.~7]{CLM16_Lorenz9D}.}\label{Fig_slowL80_vs_NN}
\end{figure}

The significant sensitivity observed in capturing the low-frequency content with nearby neural parameterizations (as measured by their loss functions) requires further investigation. Since these variations in Figure \ref{Fig_BE_vs_slowNN_closures} solely stem from how training, validation, and test sets are chosen, we conduct in this Section a statistical analysis of key features of the L80 dynamics in HLF regimes. Our focus is on the irregular lobe transitions exhibited by HLF solutions. For comparison, we also analyze lobe transitions in the slow chaotic regime of Figure \ref{Fig_sojourn}D, where neural parameterizations perform well and learn the closure effectively. Notably, Figure \ref{Fig_slowL80_vs_NN} demonstrates that for the slow chaotic regime, high-accuracy neural closures are readily achievable, with skills comparable to those obtained using the BE manifold.

To gain a deeper understanding of lobe transition statistics in the slow chaotic and HLF regimes, we performed high-resolution simulations of the L80 model for each regime. Each simulation spanned a $500$-year period, integrating the L80 dynamics with a timestep of 0.75 minutes. This corresponds roughly to an interval of size $730,000 \times T_{GW}$, where $T_{GW}$ is the dominant period of the gravity waves in the model.

In each regime, the L80 attractor exhibits two lobes.   This is  shown in the $(y_2, y_3)$-projection for the HLF regime (Figure \ref{Fig_sojourn}A) and in the $(y_1, y_3)$-projection for the slow chaos regime (Figure \ref{Fig_sojourn}D). The latter evokes the Lorenz 63 ``butterfly attractor" \cite{lorenz1963deterministic}, consistent with the L80 dynamics devoid of fast motion for this Rossby number (geostrophic motion). The former attractor, more fuzzy, exemplifies the presence of fast dynamics riding the slow, geostrophic motion.

In each case, these lobes are essentially separated by the vertical line $y_3 = 0$. Numerical integration of the L80 model reveals that the visit of the right lobe comes with $y_3(t)$ getting greater than some threshold value $y_b$, while the visit of the left lobe comes with  $y_3(t)$ getting smaller than $y_a=-y_b.$  A close inspection of the solution in the HLF  case reveals that the choice of $y_b=0.2$ constitutes a good one to identify the sojourn of the dynamics within one lobe from the other. This choice leads furthermore to an interval  $(-y_b,y_b)$ that provides a good bound of the bursts of fast oscillations crossing the vertical line $y_3 = 0$ in the $(y_2,y_3)$-plane (``gray" zone).

To count the transitions from one lobe to the other one thus proceeds as follows. Given our 500-yr long simulation of $y_3(t)$ we first find the local maxima and minima that are above $y_b$ and below $y_a$, respectively. No transition occurs between consecutive such local maxima or minima. A transition occurs only when a local maximum above $y_b$ is immediately followed by a local minimum below $y_a$ or vice versa. If a local maximum is immediately followed by a local minimum, the intermediate time instant at which the trajectory goes below zero is identified as the transition instant, and the other way around if a local minimum is immediately followed by a local maximum. These transition times characterized this way allow us to count the sojourn times in a lobe and display the distribution of these sojourn times shown in  Figures \ref{Fig_sojourn}C and \ref{Fig_sojourn}F.

These lobe sojourn time distributions reveal a striking difference between the HLF and slow chaotic regimes. 
In the HLF case, we observe indeed that the solution can stay in one lobe for a period of time that can be arbitrarily long (see solution's segment between $t=763$ and $t=893$ shown in blue in Figure \ref{Fig_sojourn}B) albeit of probability of occurrence vanishing exponentially as shown in Figure \ref{Fig_sojourn}C. As a comparison, the transitions between the attractor's lobes occur at a much more regular pace in the slow chaotic regime (see Figure \ref{Fig_sojourn}E)  in which the solutions to the L80 model are void of fast oscillations. In this case, the distribution of sojourn times drops quickly below a 60-day duration barrier (Figure \ref{Fig_sojourn}F).

These rare events, following an exponential distribution, pose a significant challenge for developing reliable slow neural closures. They introduce diversity in the temporal patterns of the time series, which contributes to the sensitivity issues observed in Figure \ref{Fig_BE_vs_slowNN_closures}.  A random training set might be skewed towards one lobe duration more than a predefined set, leading to confusion in the learning process for the neural network.

\section{The high-frequency barrier to neural closure}\label{Sec_highfreq}
 Section \ref{Sec_slow_erg} demonstrated that using neural networks to parameterize the slow dynamics of HLF solutions can lead to sensitivity issues in online prediction (Figure \ref{Fig_BE_vs_slowNN_closures}). This sensitivity arises from rare events associated with irregular lobe sojourn durations, as shown in Figures \ref{Fig_sojourn}B and \ref{Fig_sojourn}C. In this Section, we explore another challenge: the direct parameterization of the unfiltered $\x$-components of HLF solutions. These components contain a complex mixture of both slow and fast motions, posing significant difficulties for closure with neural networks.

To illustrate this point, we learn an MLP for $\x(t)$, denoted  by $\V_{\T}$, with (the unfiltered) $\y(t)$-variable of the L80 model (Eq.~\eqref{Eq_L9D}), as input, and the {\it unfiltered} $\x$-component, $\x(t)$, as output.  Note that unlike the slow NN-parameterizations above,  the parameterization $\V_{\T}$ aims at parameterizing $\x(t)$ directly as a nonlinear mapping of $\y(t)$ without conditioning on $\z(t)$ nor filtering of any sort. The corresponding closure, called a vanilla NN-closure, consists then of Eq.~\eqref{Eq_low-passNN} in which 
 $\X_{\T_2^\ast}(\y,\Z_{\T_1^\ast}(\y))$
is replaced by $\V_{\T^\ast}(\y)$,  obtained after minimization of the following $L^2$-loss function 
\bea\label{Eq_loss_vanilla}
\mathcal{L}_{\T}(\x;\y)=\sum_j \norm{\x_{t_j}-\V_{\T}(\y_{t_j})}^2,
\eea
for which the target variable $\x(t)$ is {\it unfiltered}, i.e.~containing a mixture of fast and slow oscillations. 
To address this more challenging problem we use MLPs with a larger capacity either with more neurons and/or layers. 

\begin{figure}[bth!]
\includegraphics[width=0.48\textwidth,height=.3\textwidth]{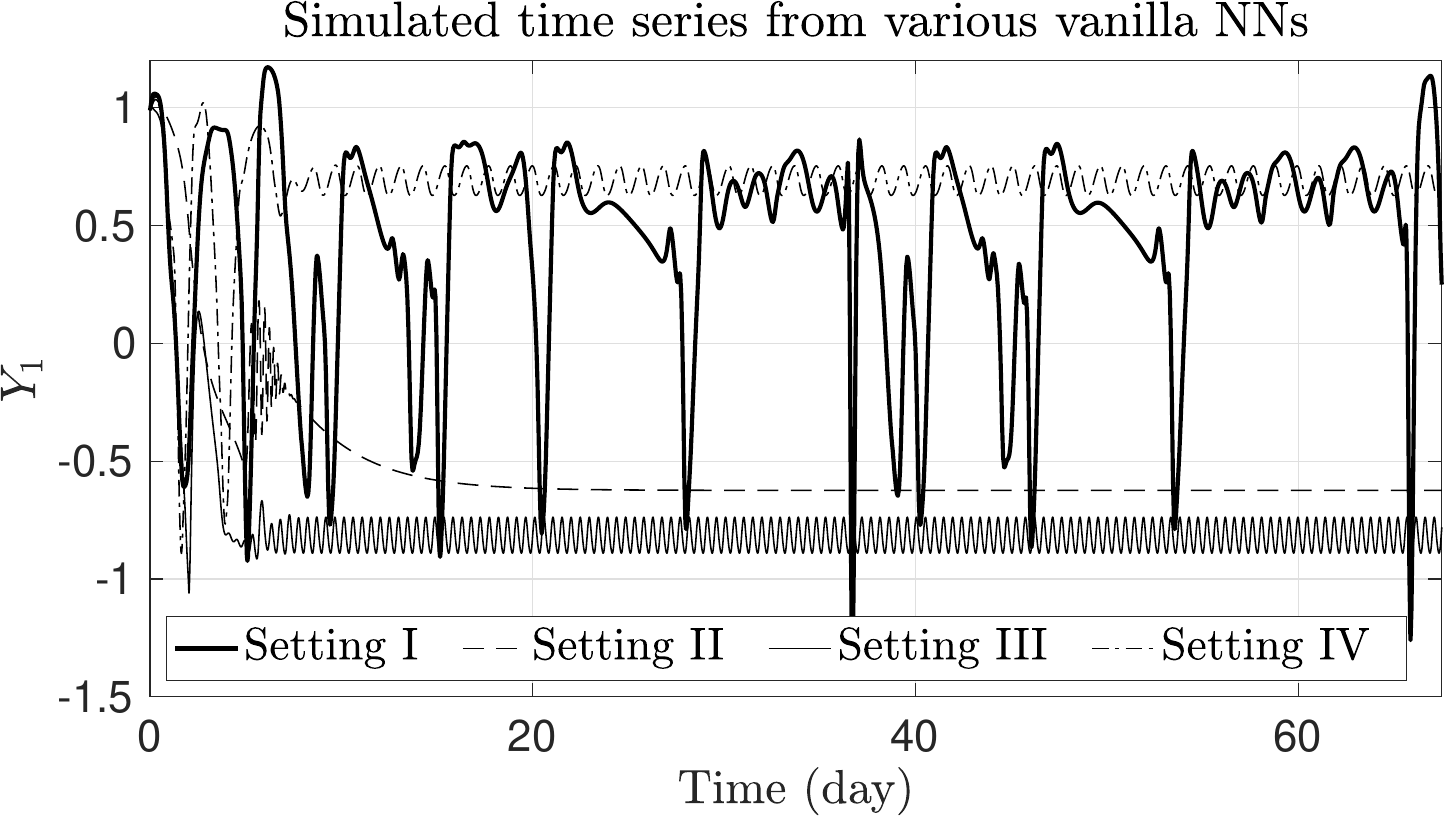}
\caption{{\bf Simulated time series from vanilla NN-closures in four different settings}. {\bf Setting I} (same as used for the results shown in Fig.~\ref{Fig_timeseries}):  one hidden layer with 20 neurons (thick solid line). {\bf Setting II}: two  hidden layers with 5 neurons in each layer (dashed line). {\bf Setting III}: two hidden layers with 10 neurons in each layer (light solid line). {\bf Setting IV}: two  hidden layers with 20 neurons in each layer (dash-dotted line). The corresponding loss function values are given in Table~\ref{tab_loss_func2}.}\label{Fig_timeseries_vanilla_NN}
\end{figure}

\begin{figure*}[bth!]
\centering
\includegraphics[width=0.95\textwidth,height=.4\textwidth]{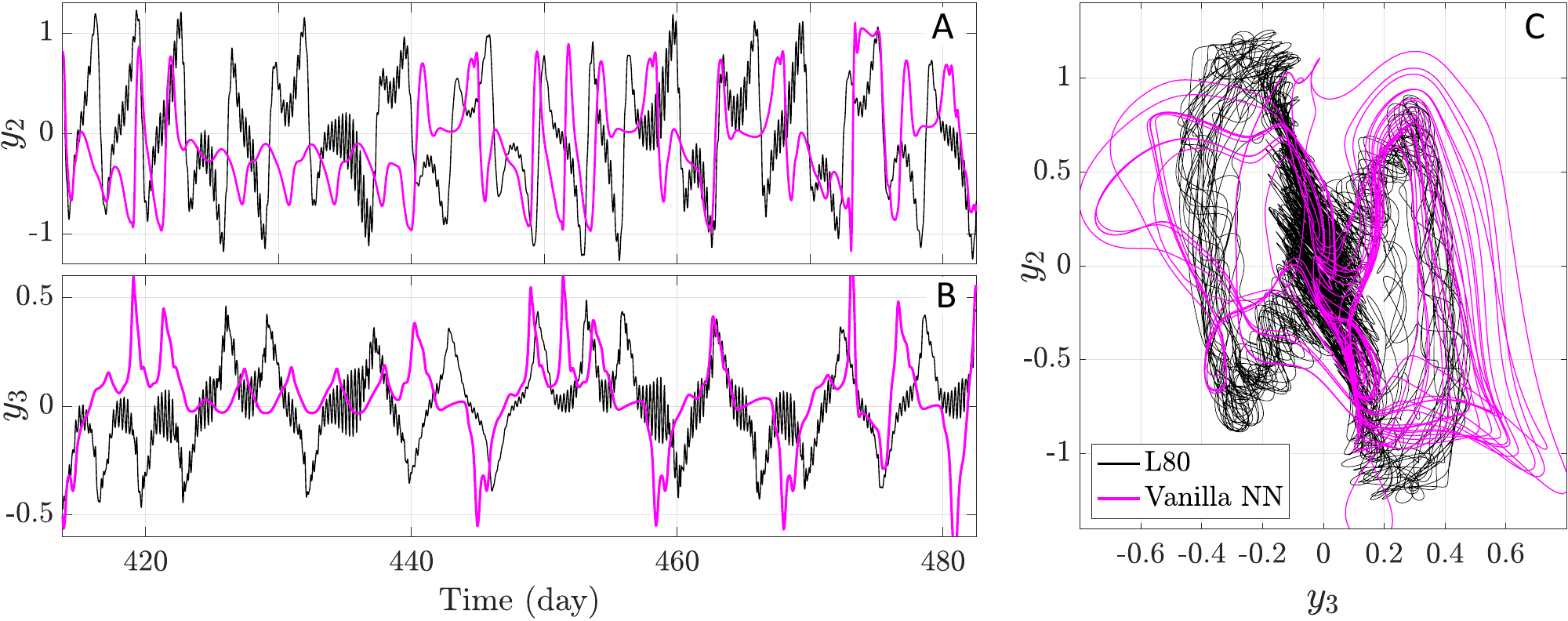}
\vspace{-2ex}
\caption{{\bf Vanilla NN-closure vs L80 dynamics}. Failure to capture the high-frequency content and symmetry of the L80 attractor. Here, is used the best performing vanilla neural network (NN1) from Setting I in Figure \ref{Fig_timeseries_vanilla_NN}}.
\vspace{-2ex}
\label{Fig_timeseries}
\end{figure*}

\begin{figure}[bth!]
\includegraphics[width=0.48\textwidth,height=.3\textwidth]{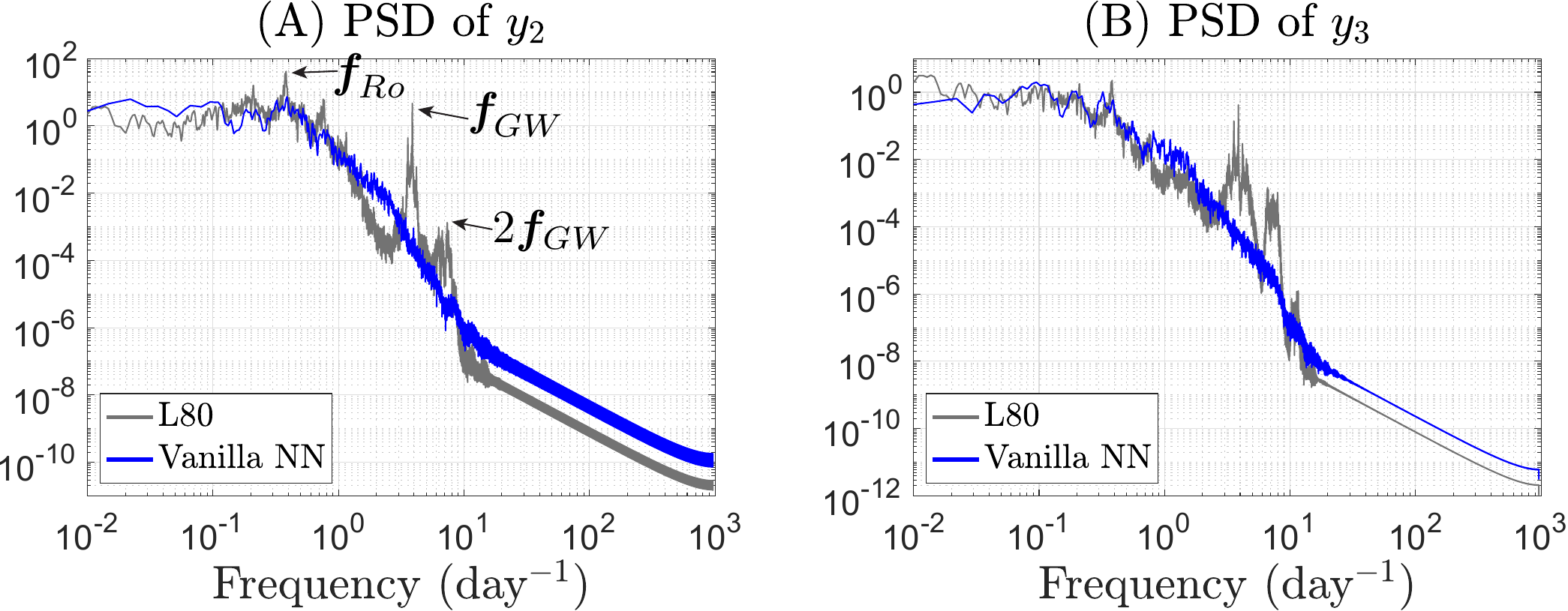} 
\caption{{\bf Power spectral density (PSD) comparison.} This figure compares the PSD of variables $y_2$ (Panel A) and $y_3$ (Panel B) for the L80 model (gray curve) and the best performing vanilla neural network closure (blue curve) from Setting I in Figure \ref{Fig_timeseries_vanilla_NN}. While the vanilla closure captures the overall spectral background of the L80 solutions well, it misses the important peaks at frequencies $f_{GW}$ and $f_{Ro}$ (and their harmonics). These frequencies correspond to inertia-gravity waves and Rossby waves, respectively.}\label{Fig_PSD0}
\end{figure}

 Interestingly, our experiments show that a neural network with just one hidden layer and $20$ neurons achieves the best closure results. Figure  \ref{Fig_timeseries_vanilla_NN} compares simulated time series from four different vanilla NN-closure settings. The setting with one hidden layer and 20 neurons partially captures the complexity of the HLF solution's temporal patterns  (Figure \ref{Fig_timeseries}A-B). However, it entirely misses the high-frequency content associated with IGWs, as evident from the power spectral density (PSD) comparison in Figure \ref{Fig_PSD0}.

\begin{figure}[bth!]
\includegraphics[width=0.48\textwidth,height=0.45\textwidth]{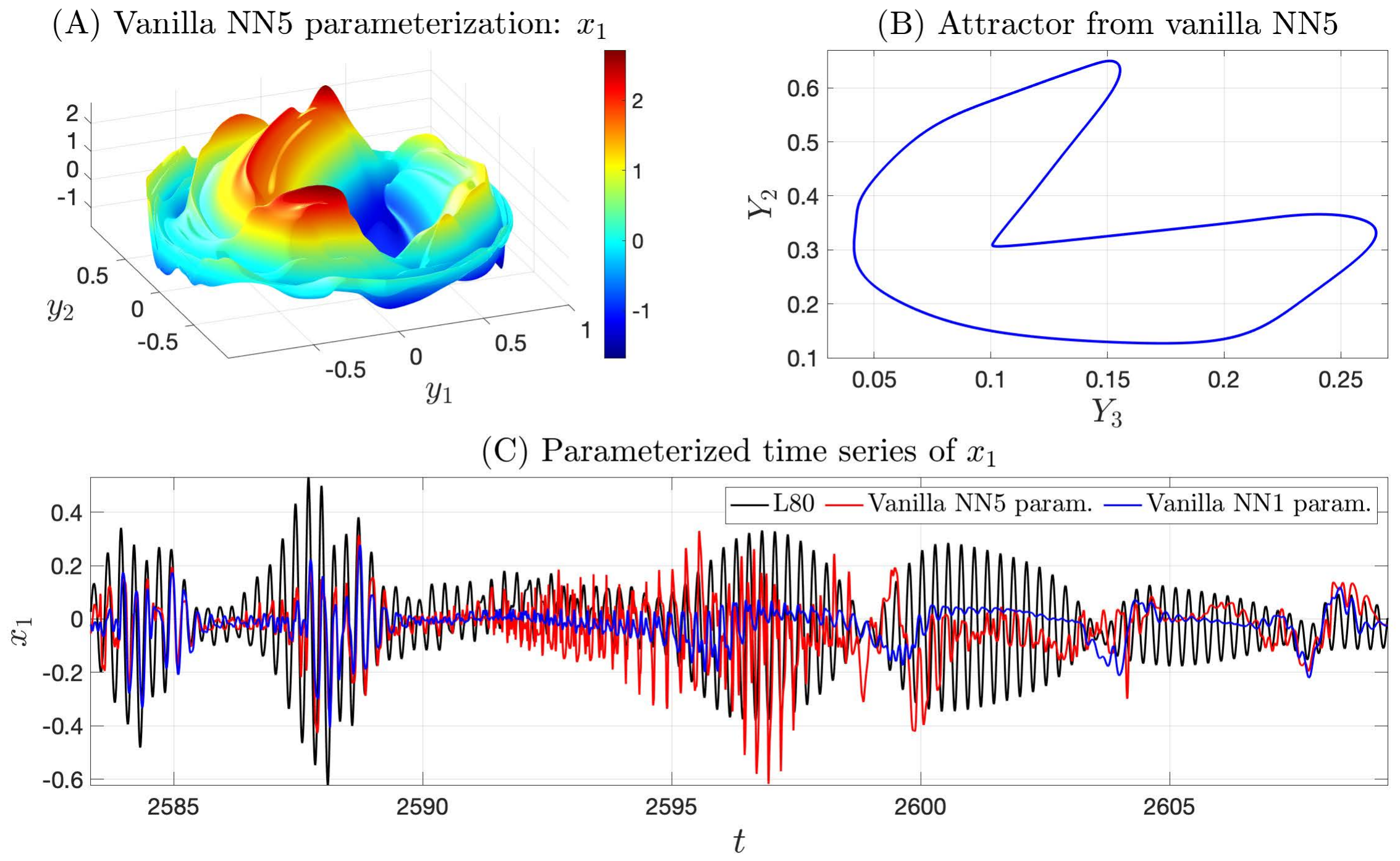} %
\caption{{\bf Panel A:} This panel shows the neural parameterization ($\V_{\T}$) with 5 layers and 20 neurons per layer (denoted as NN5) for variable  $x_1$. We use the same visualization style as Figures \ref{Fig_sensitivity}B-D. Notice the sharp gradients in the manifold, reflecting NN5's attempt to capture the high-frequency details of the HLF solutions. 
{\bf Panel B:} This panel displays the corresponding solution $(Y_2,Y_3)$ obtained using the NN5 closure.
{\bf Panel C:} Compared to the best performing vanilla neural network (NN1) from Setting I in Figure \ref{Fig_timeseries_vanilla_NN}, NN5 exaggerates the high-frequency content in the offline parameterization.}\label{Fig_BigNN}
\end{figure}

While increasing the complexity of a neural network (more hidden layers or neurons) can reduce the loss function during training, it does not guarantee better performance in the actual closure. For example, a vanilla neural network ($\V_{\T}$) with 5 hidden layers and 20 neurons per layer predicts an unrealistic, small-amplitude periodic orbit when used online in the neural closure through time-stepping (Figure \ref{Fig_BigNN}B). Additionally, it exaggerates high-frequency content in the solutions it generates offline (see Figure \ref{Fig_BigNN}C  and Table \ref{tab_loss_func2}).

Our results highlight the limitations of using a vanilla neural network closure  to directly capture the fast dynamics of the L80 system using the ``slow" variable $\y$. This approach relies on potentially complex, non-linear functions encoded by MLPs, but struggles to represent the system's multiscale dynamics accurately. This issue is similar to the spectral bias problem observed in standard neural networks for function fitting \cite{hochreiter1997long}, where they prioritize capturing low-frequency features. However, the challenge here is more complex. The goal is to learn the neglected ``fast"  variables and their high-frequency content offline,  so the online solution through the NN-closure can reproduce both the mixture of slow and fast motions of the original system. This includes capturing global geometric features like the attractor's shape and symmetry. As shown in Figure \ref{Fig_timeseries}C, vanilla NN-closures often distort these features compared to the true L80 attractor.

To address the limitations of feedforward neural networks (vanilla NNs) to close the L80 dynamics in HLF regimes, one route to explore would be to incorporate memory effects using architectures like Long-Short Term Memory (LSTM) networks \cite{hochreiter1997long}. LSTMs have demonstrably achieved model reduction in various contexts (e.g., \cite{gupta2021neural,harlim2021machine,lu2022discovering}). This success can be attributed to theoretical underpinnings from dynamical systems theory (Takens' delay embedding theorem \cite{takens1981detecting}) and statistical mechanics (Mori-Zwanzig formulation \cite{zwanzig_memory_1961,Mor65,Chorin_al02,santos2021reduced,Lucarini_Chekroun2023}). 
Additionally, we mention recent approaches combining  Takens’ 
embedding with Koopman operator theory and sparse regression
to obtain linear representations of nonlinear dynamics \cite{brunton2017chaos}.

However, as highlighted in \cite{Lucarini_Chekroun2023}, memory effects might not be crucial for achieving effcient closure of solutions in the HLF regime. Studies have shown that using  the BE manifold for capturing the geostrophic motion and a network of stochastic oscillators for IGWs can achieve high accuracy without recurrent architectures like LSTMs \cite{CLM21_BE}. This, along with the challenges of rare events discussed earlier, raises questions about whether LSTMs or other recurrent networks are necessary to reproduce the intricate multiscale dynamics of $\y$ using a closed model (like in  \cite{CLM21_BE}) built with these components (LSTMs).

\begin{table}[h] 
\caption{{\bf Loss function evaluations}. In this table are reported the loss values corresponding to the vanilla NN-closures shown in Figure \ref{Fig_timeseries_vanilla_NN}. Note that the underlying loss function is that defined in Eq.~\eqref{Eq_loss_vanilla}.}
\label{tab_loss_func2}  
\begin{tabular}{c| c c c c c c}
\toprule
Epochs & 10 & 50 & 100 & 300 & 500 \\
\hline
Setting I loss ($\times 10^{-2}$) & 2.62 & 2.54 & 2.52 & 2.49 & 2.49 \\
Setting II loss ($\times 10^{-2}$) & 2.74 & 2.67 & 2.66 & 2.64 & 2.64 \\
Setting III loss ($\times 10^{-2}$) & 2.72 & 2.45 & 2.44 & 2.43 & 2.43  \\
Setting IV loss ($\times 10^{-2}$) & 2.42 & 2.33 & 2.32 & 2.30 & 2.30 \\
\bottomrule 
\end{tabular}
\end{table}

\section{Discussion}
Our findings, particularly the interplay between rare events and the multiscale nature in HLF regimes, highlight the challenges that machine learning can face for accurate closure of  geophysical flows in which geostrophic and ageostrophic  motions  interact strongly. 
As extreme weather events and non-Gaussian statistics become more prevalent with climate change 
\cite{raveh2015large,trenberth2015attribution,swain2018increasing,galfi2021fingerprinting,seneviratne2021weather},
this study underscores that significant hurdles remain despite the recent advancements in neural parameterizations. Reliable parameterizations that robustly capture rare events are crucial. In this regard, incorporating rare event algorithms  \cite{cerou_2007,dematteis2019experimental,ragone2018computation,galfi2021applications,simonnet2021multistability} could be beneficial. By simulating rare events offline, these algorithms could improve the sampling of distribution tails, leading to better trained neural networks.

This study contributes new insights into the challenges of closing the Lorenz 80 model using data-driven methods, particularly in high Rossby number regimes ($Ro>Ro^\ast$). Compared to other Lorenz models, like the less challenging Lorenz 96 model \footnote{\url{https://raspstephan.github.io/blog/lorenz-96-is-too-easy/}}, the L80 system has received less attention for closure tasks. However, the recent stochastic closure approach by \cite{CLM21_BE} for these demanding regimes provides a valuable benchmark for future research. We hope this work encourages further exploration of the L80 model as a meaningful testbed for developing and comparing closure ideas.

\section*{Data availability statement}
The data that support the findings of this study are available upon reasonable request from the authors.

\begin{acknowledgements}
The authors sincerely thank the reviewers for their insightful comments and careful reading of our manuscript. Their feedback has been invaluable in helping us improve the clarity and focus of the article.
This work has been partially supported by  the Office of Naval Research (ONR) Multidisciplinary University Research Initiative (MURI) grant N00014-20-1-2023, by the National Science Foundation grant DMS-2108856, and by  the European Research Council (ERC) under the European Union's Horizon 2020 research and innovation program (grant agreement no. 810370).  We also acknowledge the computational resources provided by Advanced Research Computing at Virginia Tech. 
\end{acknowledgements}

\appendix

\section{HLF solutions and the slow motion learning}  \label{Sect_Appendix_HLF}
The high-low frequency (HLF) solutions used in this article are those reported in \cite[Fig.~7]{CLM21_BE}. These solutions are obtained from the parameters used in Lorenz's original paper \cite{Lorenz80} except $F_1$ chosen to be $F_1=3.027\times 10^{-1}$ as identified in \cite{CLM16_Lorenz9D}; see the Materials and Methods section in \cite{CLM21_BE} for details.

\begin{figure*}[bth!]
\includegraphics[width=0.95\textwidth,height=.45\textwidth]{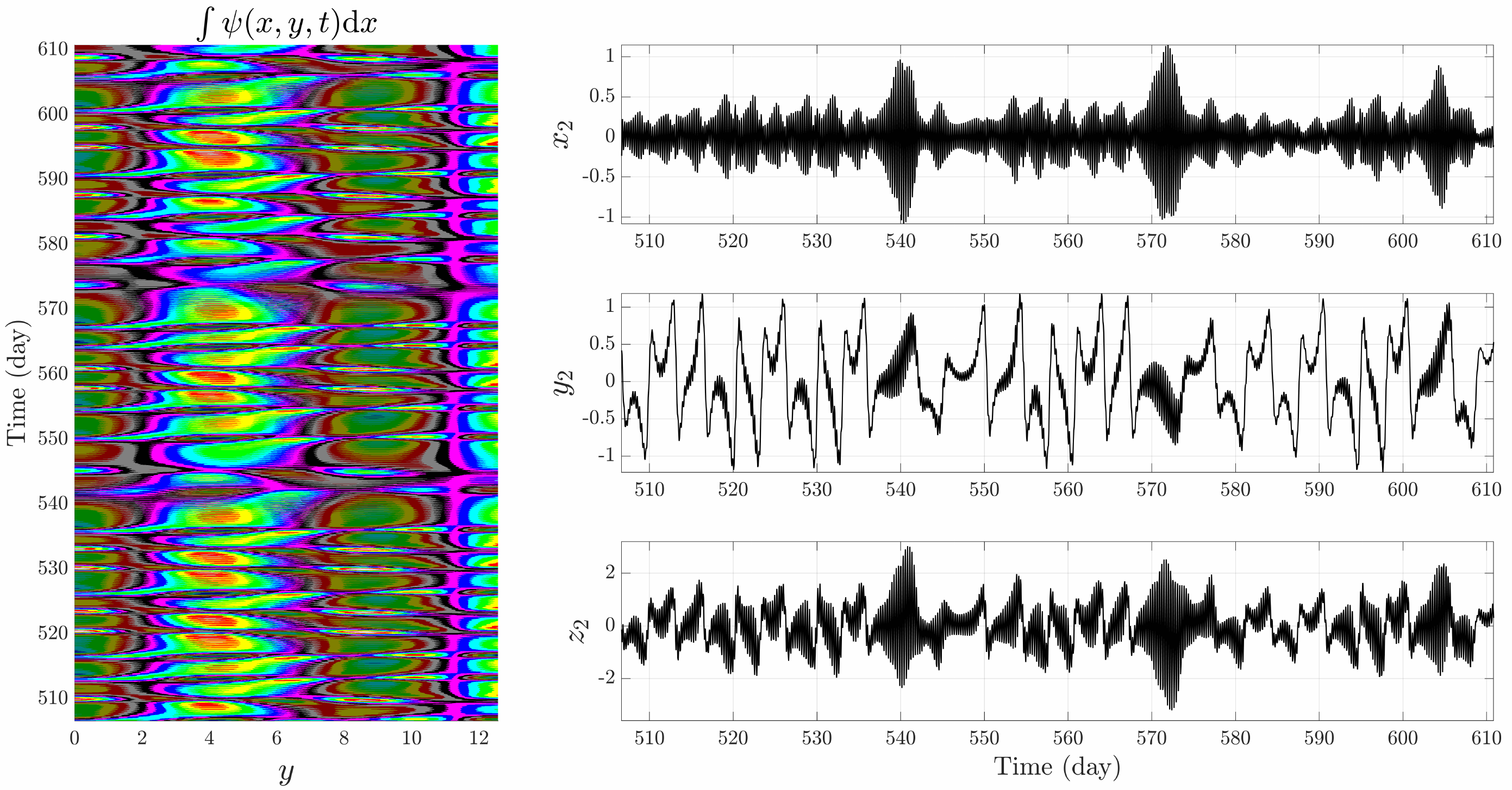} 
\caption{{\bf HLF solutions}. {\bf Left panel}: Hovm\"oller plot of the streamfunction (averaged over the $x$-direction). {\bf Right panels}: a few corresponding time series. Note the energetic bursts of  fast oscillations corresponding to spontaneous bursts of IGWs. A vanilla NN consists of seeking a feedforward neural network (FNN) mapping  
the $\y$-components of the L80 model (Eq.~\eqref{Eq_L9D}) to the $\x$-components. These fast energetic bursts are a serious barrier to learning with FNNs. The streamfunction $\psi$ is constructed from the $\y$-components of the L80 model solution according to $\psi(x,y,t) = \sum_{j=1}^3 y_j(t) \cos(\alpha_j^1 x)\cos(\alpha_j^2 y)$ where the spatial variables $x$ and $y$ (not to be confused with $\x$ and $\y$ in the L80 model) takes value in a square domain $[0, L]\times [0,L]$ with $L = 4\pi$ and the vectors $\boldsymbol{\alpha}_j = (\alpha_j^1, \alpha_j^2)$ ($j=1,2,3$) are chosen to satisfy the conditions given by \cite[Eqs.~(16)--(17)]{Lorenz80}. Following \cite{Lorenz80}, we chose $\boldsymbol{\alpha}_1 = (\sqrt{2}/2,\sqrt{2}/2)$, $\boldsymbol{\alpha}_2 = ((\sqrt{2} - \sqrt{6})/4, (\sqrt{2} + \sqrt{6})/4)$, and $\boldsymbol{\alpha}_3 = -(\boldsymbol{\alpha}_1 + \boldsymbol{\alpha}_2)$.}\label{Fig_combo_vga}
\end{figure*}

As shown in Figure \ref{Fig_combo_vga}, for this parameter regime, the HLF solutions contain a mixture of slow and fast oscillations in each variable $\x$, $\y$, and $\z$ of the L80 model that causes serious difficulties for closure \cite{CLM21_BE}. The dominant frequency of the Rossby wave content in the HLF solutions is $f_{Ro} = 0.31 \,\text{day}^{-1}$ ($T_{Ro} = 3.2 \,\text{days}$) and that of the inertia-gravity wave (IGW) content is  $f_{GW} = 3.76 \,\text{day}^{-1}$ ($T_{GW} = 6.3 \,\text{hours}$).

To learn a neural parameterization of the slow motion,  the weights and biases of the NNs are updated according to a   Levenberg-Marquardt (LM) optimization \cite{hagan1994training}. The LM algorithm is known to be efficient for small or medium-scaled problems \cite[Chap.~12]{wilamowski2018intelligent}, especially when the loss function is just a mean squared error, which is the case here. This algorithm is sufficient to obtain loss functions with small residuals; see Table \ref{tab_loss_func}.

\section{The BE manifold and BE closure} \label{Sect_Appendix_BE}   
For consistency, we recall from \cite{CLM16_Lorenz9D} the derivation of the BE manifold that serves as our parameterization baseline. 
Mathematically, the BE manifold aims at reducing the L80 model to a 3D system of ODEs, by means of nonlinear parameterization of the variables   
$\x=(x_1, x_2, x_3)^{\textrm{T}}$ and $\z=(z_1, z_2, z_3)^{\textrm{T}}$, in terms of the variable $\y=(y_1, y_2, y_3)^{\textrm{T}}$; see \cite{Gent_McWilliams82}.  
By analyzing the order of magnitudes of the different terms in  the $x_i$-equations and after rescaling following \cite{CLM16_Lorenz9D}, we arrive to the following parameterization of the $\z$-variable in terms of the rotational $\y$-variable 
\be \label{BE_z_in_y}
z_i= G_i(\y)=y_i - \frac{2c^2}{a_i} y_jy_k.
\ee

Further algebraic manipulations show that 
under an invertibility condition of a matrix $M(\y,G(\y))$ conditioned on the $\y$-variable, one obtains (implicitly) $\x$ as a function
$\Phi$ of $\y$ given by 
\be\label{Eq_phi}
\Phi(\y)=\big[M(\y,G(\y))\big]^{-1}\begin{pmatrix}
d_{1,2,3}(\y,G(\y)) \\
d_{2,3,1}(\y,G(\y))) \\
d_{3,1,2} (\y,G(\y)))
\end{pmatrix},
\ee 
where the $d_{i,j,k}$ are given explicitly; see \cite{Gent_McWilliams82,CLM16_Lorenz9D}. 
The function $\Phi (\y)=(\Phi_1(\y),\Phi_2(\y),\Phi_3(\y))^{\textrm{T}}$
corresponds to the {\it BE manifold}, it is aimed to
provide a nonlinear parameterization between $\x$ and $\y$ when the latter exists. 

The BE closure is then 
\bea\label{BE_in_y}
 \frac{\d y_i}{\d \tau} &=  -  a_i^{-1} a_kb_k \Phi_j(\y) y_k - a_i^{-1}a_j b_j y_j \Phi_k(\y) \\
 & \quad + ca_i^{-1}(a_k-a_j)y_j y_k -\Phi_i(\y)- \nu_0a_i y_i,
\eea
for which $(i,j,k)$ lies in  $\{ (1,2,3), (2,3,1), (3,1,2)\}.$

\bibliographystyle{aipnum4-1}

\end{document}